\documentclass[10pt]{article}

\usepackage[latin1]{inputenc}
\usepackage[T1]{fontenc}
\usepackage{amsfonts,amssymb,graphicx,psfrag}
\usepackage{euscript,epsfig}

\newtheorem{theo}{Theorem}
\newtheorem{prop}[theo]{Proposition}

\newtheorem{lem}{Lemma}[section]

\newenvironment{demon}[1]{{\flushleft \bf Proof #1: }}{\hfill $\square$ \vspace{5mm}}
\newenvironment{demo}{{\flushleft \bf Proof: }}{\hfill $\square$ \vspace{5mm}}
\newenvironment{dem}{{\flushleft \bf Sketch of the proof: }}{\hfill $\square$ \vspace{5mm}}

\newcommand{\R}{\mathbf R}
\newcommand{\N}{\mathbf N}
\newcommand{\M}{\widetilde{M}}
\newcommand{\D}{\partial\widetilde{M}}
\newcommand{\DD}{{\partial}^2\widetilde{M}}
\newcommand{\DDR}{{\partial}^2\widetilde{M}\times\R}
\newcommand{\Z}{\mathbf Z}

\newcommand{\un}{\mathbf 1}

\newcommand{\chap}[1]{\widehat{#1}}

\newcommand{\mcb}{\overline{\mu}_{H^+}^f}

\newcommand{\tm}{T^1M}
\newcommand{\ttm}{T^1\widetilde{M}}
\newcommand{\w}{\widetilde{W}}
\newcommand{\ffi}{\widetilde{\Phi}}

\begin{document}

\title{On quasi-invariant transverse measures for the 
horospherical foliation of a negatively curved manifold}


\author{ Barbara SCHAPIRA\\
\\
\em  MAPMO,
Universit\'e d'Orl\'eans,
Rue de Chartres,\\
\em BP 6759,
45067 Orl\'eans cedex 2, France\\
\em schapira@labomath.univ-orleans.fr\\}

\maketitle

\begin{abstract}
\noindent
If $M$ is a compact or convex-cocompact negatively curved manifold, we associate to
any Gibbs measure on $\tm$ a quasi-invariant transverse measure for the horospherical foliation, and prove that this measure is uniquely determined 
by its Radon-Nikodym cocycle. (This extends the Bowen-Marcus unique ergodicity result for this foliation.) 
We shall also prove equidistribution properties for the leaves of the foliation w.r.t. these Gibbs measures.
We use these results in the study of invaiant meausres for horospherical foliations on regular covers of $M$.\\
\end{abstract}

\vfill

\noindent
{\bf Primary AMS classification:} 37D40, 37C85, 37A20, 22F05 \\

\noindent
{\bf Secondary  AMS classification:} 37D35, 32Q05, 58Hxx,\\

\vspace{0.5cm}


\vfill
\newpage

\section{Introduction and statement of results}

Let $M= \Gamma\backslash \widetilde{M}$  be a complete riemannian manifold with
 pinched negative curvature, 
with universal cover $\widetilde M$ and fundamental group $\Gamma$. Then
 the  unit tangent bundle  $\ttm$ of $\widetilde M$ carries the {\em horospherical
foliation} $\widetilde{\cal W}^{su}$, whose leaves are  the unit
vectors normal to horospheres and pointing outward. Passing to the quotient
leads to a
foliation ${\cal W}^{su}$ of $\tm$ which is
 the  strong unstable foliation
of the geodesic flow of $M$. In this paper, we shall construct
quasi-invariant measures for
this foliation, and study their ergodic properties.

\noindent
Recall that for a compact manifold $M$,  it was  proved by Bowen and Marcus
that the foliation ${\cal W}^{su}$ is uniquely ergodic:
there exists, up to a multiplicative constant, a unique  transverse measure
for the
foliation which is invariant under holonomy.
Their proof used symbolic dynamics and showed that the transverse measure
is induced by
the  Bowen-Margulis measure $m^0$ of $T^1M$, which is also the measure of
maximal entropy of
the geodesic flow \cite{BM}.

\noindent
In the last decade, it was realized that the  measure $m^0$ could be
constructed  in a
more geometrical way
  using the  Patterson measure on the boundary at
infinity of
$\widetilde{M}$
\cite{kaim2} \cite{kaim}. This   construction allows to go beyond compact
manifolds,
and,   for  instance,   Roblin was able to give a purely geometrical proof
of the
unique ergodicity of ${\cal W}^{su}$  for a convex-cocompact manifold,
under the assumption
that the geodesic flow is topologically mixing \cite{roblin}.

\noindent
In a similar way, one can associate to a   H\"older function $f : T^1M\to \R$ a
measure $m^f$ on $T^1M$, which is the equilibrium state of $f$ when $M$ is
compact or
convex-cocompact,  the
Bowen-Margulis measure
$m^0$ corresponding to $f\equiv 0$ see e.g.
\cite{ham}\cite{Ledrappier}\cite{coudene}.  We shall first  see that this
measure
$m^f$ induces a  {\em quasi-invariant  transverse measure} for the foliation
${\cal W}^{su}$.  We  use here the definition of a quasi-invariant measure of a
foliation  as introduced by Connes in \cite{Connes}~ in relation with the
theory of
$C^*$-algebras associated to a foliation.  First, we recall that a {\em cocycle}
for a
foliation is a map $\rho$ defined on the set of pairs of points in the same
leaf
such that $\rho(u,v)+\rho(v,w)=\rho(u,w)$, and that a {\em transverse measure}
$\mu =\{\mu_T\}$
associates to any tranversal $T$ to the foliation a Borel measure $\mu_T$
supported on $T$.
Then $\mu$ is said to be  quasi-invariant  if there exists a cocycle $\rho$
such
that for each holonomy map $\zeta$ between two transversals $T$ and $T'$,
$\zeta_*\mu_T$ is
 absolutely continuous with respect to $\mu_T'$ with   Radon-Nikodym
derivative given by
 $$
\frac{d\zeta_*\mu_T}{d\mu_{T'}}(\zeta x) =\exp(\rho(x, \zeta x)).$$
Note that if $\mu'=\{\mu'_T\}$ is {\em equivalent} to $\{\mu_T\}$ (i.e. there
exists a Borel  map
$\psi:\tm\to\R$ such that, for all transversal $T$,
$d\mu_T(v)=\exp\psi(v)\,d\mu'_T(v)$),
then the Radon-Nikodym cocycle $\rho'$ of $\mu'$ is {\em cohomologous} to
$\rho$
(there exists a Borel  map $R:\tm\to\R$, such
that $\rho(v,w)=\rho'(v,w)+R(v)-R(w)$).\\
Let us call two H\"older functions
$f$ and $f'$ {\it equivalent} if there exists a H\"older map
$\phi:\tm\to\R$, differentiable in the direction of the geodesic flow, and a
constant $c$
such that $f=f'+c+X.\phi$, with $X:TM\to TTM$ the geodesic spray.\\
Our starting point is the following Proposition: it provides a large family
of
cocycles which are
Radon-Nikodym cocycles of quasi-invariant measures.

\begin{prop}\label{prems}To each H\"older map $f:\tm\to\R$ is associated
an explicit H\"older cocycle $\rho^f$ for the foliation
${\cal W}^{su}$,
and a transverse measure $\mu^f=\{\mu_T^f\}$,
which is quasi-invariant with cocycle $\rho^f$.
Moreover, the cohomology class of $\rho^f$ and the equivalence class of
$\mu^f$
depend only the equivalence class of $[f]$.
\end{prop}

\noindent
For $f\equiv0$, the cocycle   $\rho^0$ is trivial, and one recovers a
transverse invariant
measure $\mu^0$. For compact or convex-cocompact manifolds, one knows that this measure is
the unique
invariant mesure of the foliation, and thus we are led to the following question: is
the transverse measure $\mu^f$ the unique quasi-invariant  measure with the
given
cocycle
$\rho^f$?

\noindent
In the sequel of this work, we will always make the following assumptions:\\

\noindent
{\bf Assumptions: }{\em The fundamental group  $\Gamma=\pi_1(M)$ is
cocompact or
convex-cocompact,
and in the second case, the geodesic flow is topologically mixing on its
nonwandering set.} \\

\noindent
Under these assumptions, our main result is:

\begin{theo}\label{1bis} Let $[\rho^f]$ be a cohomology class of cocycles
associated to a class $[f]$ of H\"older maps.
Then for all $\rho\in[\rho^f]$, there exists, up to a multiplicative constant, a
unique
transverse measure $\mu$ quasi-invariant with Radon-Nikodym cocycle
$\rho$, and it is equivalent to
$\mu^f$.
\end{theo}

\noindent
When $f\equiv 0$  we recover the unique ergodicity of the foliation.
Theorem \ref{1bis} was already proved by Babillot-Ledrappier 
\cite{mbab} for a compact manifold using symbolic dynamics,
but we give here a geometrical proof, which allows
to prove further results as the equidistribution property
of leaves of the foliation (see below Theorem \ref{convergence}).\\

\noindent
Theorem \ref{1bis} can be reinterpreted in terms of the
action of the fundamental group $\Gamma$ of $M$ on the boundary at infinity
$\D$
of $\M$. Indeed the lifted
 foliation
$\widetilde{\cal W}^{su}$ of $\ttm$
 admits a nice set of leaves: the space ${\cal H}$ of horospheres,
which can be
parametrized as $\D\times\R$.
An invariant transverse measure for ${\cal W}^{su}$ lifts to an invariant
transverse measure for $\widetilde{\cal W}^{su}$ which is
$\Gamma$-invariant, and thus   induces in a canonical way a $\Gamma$-invariant
measure on ${\cal H}$. In the case of quasi-invariant
measures, such a correspondance can be made by chosing a global tranversal
of $\widetilde{\cal
W}^{su}$, and thus depends on this choice.
Our construction is slightly different : we associate to each
quasi-invariant transverse measure $\{\mu_T^f\}$ a quasi-invariant measure on
${\cal H}$ with the additional property that its
Radon-Nikodym  cocycle
$$c^f : \Gamma\times {\cal H}\mapsto\R$$ is in fact a cocycle for the action of
$\Gamma$ on $\D$:
for all $(\xi,s)\in\D\times\R\simeq {\cal H}$, and $\gamma\in\Gamma$,
$$c^f(\gamma,(\xi,s))=c^f(\gamma,\xi).$$
To avoid confusions, we will denote by $[c^f]_{\cal H}$ and $[c^f]_{\D}$ the
respective
cohomology classes of $c^f$ as a cocycle for the action of $\Gamma$ on ${\cal
H}$ and on $\D$.

\noindent
A result of Ledrappier, \cite{Ledrappier},
shows that any H\"older cocycle $c$ on $\D$ can be written as $c=c^f$, for a
convenient H\"older map $f:\tm\to\R$.
Thus we can reformulate Theorem \ref{1bis} in the following way:\\

\noindent
{\bf Theorem 2 bis} {\it Let  $[c]_{\cal H}$ be the cohomology class
on ${\cal H}$ of a H\"older cocycle $c$ on $\D$. There exists, up to a
multiplicative
constant, a unique
class $[\chap{\mu}]$ of $\Gamma$-quasi-invariant measures on ${\cal H}$
with respect to this class $[c]_{\cal H}$, and they are ergodic.
}\\

\noindent
This result  leads to the question whether
there exist cohomology classes of cocycles on ${\cal H}$ which do not contain
cocycles arising from cocycles on $\D$,
and if the $\Gamma$-quasi-invariant measures on ${\cal H}$ with such a cocycle
can be characterized. \\

\noindent
In the case of a flow, the property of unique ergodicity
implies (and is even equivalent to)
the equidistribution of all orbits to the unique invariant measure.
In our set-up,  we have no flow, so we have to
consider more general means on leaves of the foliation.
Then we shall see that the above implication is still true.

\begin{theo}\label{convergence} Let $f:\tm\to\R$ be a  H\"older map. There
exists on
each strong unstable horosphere $H^+$ a  measure $\overline{\mu}^f_{H^+}$,
and a distance
$d_{H^+}$, such that
for all non wandering vector
$u\in\tm$, the mean $M_{r,u}$ on the ball $B^+(u,r)$ w.r.t.
the measure $\overline{\mu}^f_{H^+}$
converges weakly to the equilibrium state  $m^f$ associated to $f$ when
$r\to\infty$.

\end{theo}

\noindent
We will also see that the equilibrium state $m^f$ is locally equivalent to the
product $\chap{\mu}^f\times\overline{\mu}^f_{H^+}$.

\

\noindent
In fact, we shall deduce this Theorem as a corollary of the more general
Theorem \ref{p7} of equidistribution of sequences of sets $(E_n)_{n\in\N}$
on leaves which satisfy a certain  growth property.\\

\noindent
As another application of the preceding results, we shall address the
problem of determining the $\bar{\Gamma}$-invariant measures on ${\cal H}$,
with $\bar{\Gamma}\triangleleft\Gamma$ a normal subgroup of $\Gamma$.
Babillot and Ledrappier proved in \cite{mbab} that when $\Gamma$ is cocompact,
and $\Gamma/\bar{\Gamma}$ is isomorphic to $\Z^d$, each
cohomology class $[\alpha]$ of $1$-forms on $M$, which vanishes on loops of
$\bar{\Gamma}$,
induces a $\bar{\Gamma}$-invariant measure on ${\cal H}$, and that these
measures are $\bar{\Gamma}$-conservative and ergodic.
However, it is still an open question to see whether
the measures they constructed are the only $\bar{\Gamma}$-invariant and
ergodic
measures on ${\cal H}$.
Here is a partial answer under an additional assumption: for a
general normal subgroup $\bar{\Gamma}$ of a cocompact or convex-cocompact
group
$\Gamma$, these measures are the only possible
$\Gamma$-quasi-invariant, $\bar{\Gamma}$-invariant, and $\bar{\Gamma}$-ergodic
measures on ${\cal H}$.

\begin{theo}\label{cover}
Let $\overline{\Gamma}\triangleleft\Gamma$ be a normal subgroup. If
$\alpha:TM\to\R$ is a closed $1$-form
vanishing on the image of $\bar{\Gamma}$ in $H_1(M,\R)$,
then the measure $\chap{\mu}^{\alpha}$ associated to the H\"older map
$\alpha_{|\tm}$ is $\bar{\Gamma}$-invariant.
Conversely, any $\overline{\Gamma}$-invariant and  $\overline{\Gamma}$-ergodic
measure $\chap{\nu}$  on ${\cal H}$, which
 is quasi-invariant by $\Gamma$, is necessarily a measure of the form
$\chap{\nu}= \chap{\mu}^{\alpha}$,
with $\alpha$ a closed $1$-form (vanishing on $\overline{\Gamma}$).
\end{theo}

\noindent
Note that in the general case $\overline{\Gamma}\triangleleft\Gamma$, it is
not
clear whether
the measures $\chap{\mu}^{\alpha}$ are $\overline{\Gamma}$-ergodic;
 see \cite{kaim3}, \cite{kaim4} for cases where it is known and other
references.\\

\noindent
The organization of the text is the following: in a first part, we introduce
some notations and all the measures that we will study.
In part \ref{equicontinuite}, we prove
the  unicity results (Theorems \ref{1bis} and \ref{1bis} bis). We deduce
then of
these results an equidistribution result for balls and for more general
sets on horospheres (Theorem \ref{p7} and Theorem \ref{convergence}).
Finally, we prove Theorem \ref{cover} in part \ref{nonsym}.  \\

\noindent
It is a pleasure to thank Vadim Kaimanovich for his very careful reading, and many useful suggestions on a preliminary version of this work.
Also, I would like to thank my supervisor Martine Babillot for her continuous and valuable 
guidance during the elaboration of this article, and Jean-Pierre Otal for several discussions.


\section{Notations, preliminaries}


\subsection{Geometry}

Let $M=\Gamma\backslash\M$ be a complete riemannian manifold with pinched negative sectional curvature, 
$\M$ its universal cover, and $\Gamma$ its fundamental group. Then $\M$ can be compactified into $\overline{M}=\M\cup\D$, where $\D$ 
is the {\it boundary at infinity} of $\M$, 
i.e the set of equivalence classes of geodesic rays
which stay at bounded distance one another.
The group $\Gamma$ acts on $\M$ by isometries, and
on $\D$ by homeomorphisms.
The {\it limit set}  $\Lambda$ of $\Gamma$ is the set of accumulation points of $\Gamma$ in $\D$: 
$\Lambda=\overline{\Gamma x}\backslash\Gamma x$, for any $x\in\M$.
We denote by $T^1M$ (resp. $\ttm$) the unit tangent bundle of $M$ (resp. $\M$), and by $\pi:\tm\to M$ (resp. $\ttm\to\M$) 
the canonical projection.
We will use the distance $d$ on $M$ (and $\M$) induced by the riemannian structure.

\noindent
The geodesic flow $\Phi=(\Phi^t)_{t\in\R}$ (resp. $\ffi$) 
associates to a pair $(t,v)\in\R\times\tm$ (resp. $\R\times\ttm$)
the tangent vector
$\Phi^t v=\dot{c}_v(t)$ at time $t$ to the unique geodesic $c_v$ of $M$ such that $\dot{c}_v(0)=v$.

\noindent
By a theorem of Eberlein \cite{eberlein}, the {\it nonwandering set} $\Omega\subset\tm$ of the geodesic flow $\Phi$ is the set of vectors $v$  
such that any lift $\widetilde{v}$ to $\ttm$ has both endpoints $c_v(\pm\infty)$ in the limit set $\Lambda$.

\noindent
The group $\Gamma$ is {\it cocompact} if $M=\Gamma\backslash\M$ is compact, and then  
 $\Lambda=\D$, and $\Omega=\tm$. 
It is {\it convex-cocompact} when $\Omega$ is compact. 
In this work, we assume that $\Gamma$ is cocompact or convex-cocompact, 
so $\Omega$ will always  be a compact set. We assume also that the geodesic flow is {\it topologically mixing} on its non-wandering set $\Omega$, 
which is equivalent to the nonarithmeticity of the length spectrum of $M$. It is known to be true 
in the case of a surface, in the constant curvature case, and in some other cases (see Dal'bo, \cite{dalbo}).

\noindent
A H\"older {\it cocycle} on $\D\times\M\times\M$ is a H\"older map $c:\D\times\M\times\M\to\R$ such that for all $\xi\in\D$, and $x,y,z$ in $\M$, 
$c_{\xi}(x,y)+c_{\xi}(y,z)=c_{\xi}(x,z)$. It is said to be $\Gamma$-invariant if it is invariant under the diagonal action of $\Gamma$ 
on $\D\times\M\times\M$.
The {\it Busemann cocycle} is defined on $\D\times\M\times\M$ by:
$$
\beta_{\xi}(x,y)=\lim_{z\to\xi}d(x,z)-d(y,z)="d(x,\xi)-d(y,\xi)".
$$
It is a continuous and $\Gamma$-invariant cocycle on $\D\times\M\times\M$. 
By abuse of notation, if $v,w$ are vectors on 
$\ttm$ with base points $x$ and $y$ on $\M$, we will often note $\beta_{\xi}(v,w)$ in place of $\beta_{\xi}(x,y)$. 

\noindent
The space of oriented geodesics can be identified to the {\it double boundary}, i.e. the set of pairs of distinct points 
$(\eta_1,\eta_2) \in  \DD:=\D\times\D\backslash \{(\xi,\xi),\xi\in\D\}$. 
The Busemann cocycle allows to give coordinates on $\ttm$ in terms of this double boundary.
More precisely, if a point $o\in\M$ is fixed, the map defined by:  
$$
v\in\ttm\to (v^-,v^+,\beta_{v^-}(v,o))\in\DDR
$$
is an homeomorphism.

\begin{figure}[!ht]
\begin{center}
\input{coordo_bord.pstex_t}
\caption{Coordinates on $\ttm$}
\end{center}
\end{figure}
\noindent
On $\DDR$, the actions of the geodesic flow and the group $\Gamma$ commute and can be written in the following way: for
any $\gamma\in\Gamma$,
$$
\gamma(v^-,v^+,s)=(\gamma v^-,\gamma v^+, s+\beta_{v^-}(o,\gamma^{-1} o))\quad\mbox{and} \quad\ffi^t(v^-,v^+,s)=(v^-,v^+,s+t).
$$
\noindent
Thus, the unit tangent bundle $\tm$ identifies to the quotient $\Gamma\backslash(\DDR)$, and the non-wandering set $\Omega$ to 
$\Gamma\backslash(\Lambda^2\times\R)$.

\noindent
A {\it  horosphere} $H\subset\M$ centered at $\xi$ is a level set of a function  $x\to\beta_{\xi}(x,y)$. 
These horospheres lift to $\ttm$: 
if $u\in\ttm$, and $H$ is the horosphere centered at $u^-$ and containing the base point  $\pi(u)\in \M$ of $u$, 
the {\it strong unstable horosphere of $u$}, 
denoted by $H^+(u)$ is the set of vectors $v\in\ttm$ such that $v^-=u^-$ and 
$\pi(v)\in H$.
It is also the set of vectors of  $\ttm$ with base points in $H$, orthogonal to $H$ 
and pointing outwards. 
Similarly, we denote by $H^-(u)$ the strong stable horosphere of $u$. We denote by ${\cal H}$ the {\it space of horospheres}.

\begin{figure}[!ht]
\begin{center}
\input{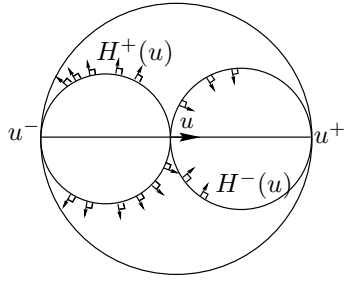}
\caption{Horospheres}
\end{center}
\end{figure}

\noindent
The homeomorphism  $\ttm\simeq\DDR$ allows to identify a horosphere $H^+ (u)$ with $\{u^-\}\times\D\backslash\{u^-\}\times\{s(u)\}$,
and the space ${\cal H}$ of horospheres with $\D\times\R$. The group $\Gamma$
acts on ${\cal H}$ by $\displaystyle \gamma(\xi,s)=(\gamma\xi,s+\beta_{\xi}(o,\gamma^{-1}o)$.
Besides, 
the horospheres $H^+(u)$ are exactly the strong unstable manifolds of the geodesic flow  $\ffi$: 
$$
H^+(u)=\w^{su}(u):=\{w\in\ttm, \, \lim_{t\to -\infty}D(\ffi^t u,\ffi^t w)=0\}.
$$ 
\noindent
Here, $D$ is the distance on $\ttm$ induced by the Sasaki metric on $T\M$.
Similarly, the strong stable horospheres  $H^-(u)$ equal the strong stable manifolds $\w^{ss}(u)$.
These strong unstable horospheres form a foliation $\widetilde{\cal W}^{su}$ of $\ttm$, 
called the {\it horospherical foliation}, 
or strong unstable foliation. 
A natural family of transversals to $\widetilde{\cal W}^{su}$ is the family of weak stable manifolds: 
$$
\w^s(v)=\{w\in\tm, \exists C=C_w>0, \forall t\geq 0, D(\Phi^t v,\Phi^t w)\leq C\}.
$$

\noindent
Viewed on $\DDR$, a transversal $\w^s(v)$ equals $\D\backslash\{v^+\}\times\{v^+\}\times\R$.
Recall that, given two (small) transversals $T$, $T'$ to the foliation included in the same chart of the foliation,
a {\it holonomy map} is a homeomorphism
from an open subset of $T$  to $T'$, which preserves the leaves of the foliation.
Note that there are natural (global) holonomy maps between two such transversals, given by:
$$
\begin{array}{rcl}
\w^s(v) & \!\!\! \stackrel{\zeta}{\longrightarrow} \!\!\! & \w^s(w)\\
u=(u^-,v^+,s(u)) & \!\!\! \mapsto \!\!\! & \zeta(u)=(u^-,w^+,s(u)).
\end{array}
$$

\noindent
The foliation $\widetilde{\cal W}^{su}$ 
induces on the quotient $\M/\Gamma$ the foliation ${\cal W}^{su}$ whose leaves are the strong unstable 
manifolds $W^{su}$ for $\Phi$. 
For the latter foliation, the holonomy maps are not defined globally. \\

\noindent
Finally, let us introduce useful families of distances.
First, the family $(d_x)_{x\in\M}$ of Gromov distances on the boundary is defined by: 
$$
\forall (\xi,\eta)\in\DD,\,d_x(\xi,\eta)=\exp\left(-\frac{1}{2}\beta_{\xi}(x,y)-\frac{1}{2}\beta_{\eta}(x,y)\right),
$$
with $y$ any point on the geodesic $(\xi,\eta)$. 
(Rigorously, they are known to be true distances if the curvature of  $M$ is less than or equal to $-1$ \cite{bourdon}; 
otherwise, one has to take convenient powers 
of the above quantities to obtain distances.) 
Using these distances, Hamenst\"adt \cite{ham2} defined on each horosphere $H^+$ a distance $d_{H^+}$ by the formula:
$$
\forall (u,v)\in (H^+)^2,\,\,d_{H^+}(u,v)
=\exp\left(\frac{1}{2}\beta_{u^+}(x,u)+\frac{1}{2}\beta_{v^+}(x,v)\right)\,d_x(u^+,v^+)\,,
$$

\noindent
where $x\in\M$ is any point (one easily checks that the above formula does not depend on $x$).
Remark that these distances can also be defined on the strong stable horospheres $H^-(v)$.
The following picture shows what they represent geometrically:
if $u,v\in H^+$, 
$d_{H^+}(u,v)=\exp(\pm d(a_u,a_v)/2)$ where $a_u$ and $a_v$ are the respective intersection points 
of $H^-( u)$ and $H^- (v)$ with the geodesic $(u^+,v^+)$ (the sign depends on the order of  $a_u$ and $a_v$ on $(u^+,v^+)$). 
\begin{figure}[!ht]
\begin{center}
\input{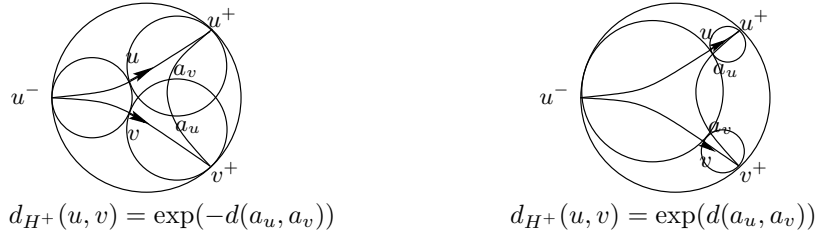}
\caption{Horospherical distance}\label{dist}
\end{center}
\end{figure}
\noindent

\noindent
Two important properties of this family are that they are
 invariant by isometries:
$$\forall\gamma\in\Gamma,\,d_{\gamma H^+}(\gamma u,\gamma v)=d_{H^+}(u,v),$$
and that they  ``explode'' when pushed by the flow: 
$$\Phi^t B^+(u,r)=B^+(\Phi^t u, re^t)$$


\subsection{Cocycles and measures associated to a H\"older map}\label{measures}

A {\it H\"older potential} is a map $f:\tm \to\R$ which is H\"older 
w.r.t. the Sasaki metric $D$ on $\tm$. 
For simplicity, we shall consider for the moment only symmetric potentials:
 $f(v)=f(-v)$ (see section \ref{chech} for the general case), 
and we denote also by $f$ its $\Gamma$-invariant lift on $\ttm$. 
If $x$ and $y$ are points in $\M$, 
we denote by $\int_x^y f$ the integral of $f$ on the unique (oriented) geodesic  
from $x$ to $y$, viewed on $\ttm$. 
If $x,y\in\M$ and $\xi\in\D$, we define the following cocycle:
$$
\rho^f_{\xi}(x,y)=\lim_{t\to +\infty}\int_x^{c(t)}f-\int_y^{c(t)}f\,=\,"\int_x^{\xi}f-\int_y^{\xi}f",
$$ 
where $(c(t))_{t\geq 0}$ is the geodesic ray $[x,\xi)$.
This function is well-defined, since when $f$ is H\"older, 
the difference between the two integrals converges. Moreover, it is $\Gamma$-invariant. 
It is a generalisation of the Busemann cocycle, since when
$f\equiv 1$, $\rho^1_{\xi}(x,y)=\beta_{\xi}(x,y)$.
Moreover, it induces a cocycle on the strong unstable foliation $\widetilde{\cal W}^{su}$: 
if $v$ and $w$ are two vectors on the same strong unstable horosphere $H^+(v)$, this cocycle is defined by:
$\displaystyle \rho^f(v,w)=\rho^f_{v^-}(v,w).$\\
Note that when $f\equiv 1$, $\rho^f\equiv 0$. By $\Gamma$-invariance, one may also define 
this cocycle for the strong unstable foliation ${\cal W}^{su}$ on $\tm$;  
indeed, if $v$ and $w$ are on the same leaf of $\tm$, 
the quantity $\rho^f(\widetilde{v},\widetilde{w})$ is $\Gamma$-invariant, so it 
does not depend on the lifts $\widetilde{v},\widetilde{w}$ chosen on the same leaf of $\widetilde{\cal W}^{su}$. 
We can also define it directly on $\tm$ by
$$\rho^f(v,w)=\lim_{t\to -\infty}\int_{v}^{\Phi_t v}f-\int_{w}^{\Phi_t w}f.$$

\noindent
Finally, if we choose an origin $o\in\M$, we can associate to $f$ the following  H\"older cocycle on $\D$: 
$$
c_o^f(\gamma,\xi)=-\rho^f_{\xi}(o,\gamma^{-1}o). 
$$

\begin{lem}\label{classescohom}The cohomology class $[c_o^f]_{\D}$ depends on $f\in[f]$, but on ${\cal H}$, the class
$[c_o^f]_{\cal H}$ depends only on $[f]$, and neither on $o$, nor on $f\in[f]$.
\end{lem}
The proof is an easy computation.
We will see after the following proposition that there is another interesting cocycle on $\D$ associated to $f$. 

\noindent
To each H\"older potential $f$ on $\tm$ is associated a unique equilibrium state $m^f$, 
that is a probability measure on $\Omega\subset \tm$ satisfying the variational principle 
(see \cite{Ledrappier}, \cite{coudene} for the construction, and \cite{kaim2}, \cite{ham} to see that it is an 
equilibrium state).
This measure is induced on the quotient $\tm=\Gamma\setminus\ttm$  by a measure $\widetilde{m}^f$ on $\ttm$. 
We will recall the two main steps of its construction. 

\begin{prop}[Ledrappier, \cite{Ledrappier}] There exists a (unique) probability measure $\nu_o^f$ on $\Lambda\subset\D$, 
which has no atoms, and is ergodic and $\Gamma$-quasi-invariant with respect to the following cocycle:
$$(\gamma,\xi)\to \delta^f\beta_{\xi}(o,\gamma^{-1}o)-\rho_{\xi}^f(o,\gamma^{-1}o) .$$
The constant $\delta^f$ is the topological pressure of $f$.
Moreover the equivalence class $[\nu_o^f]$ of $\nu_o^f$ and the cohomology class on $\D$ of its Radon-Nikodym cocycle 
$\delta^f\beta_{\xi}(o,\gamma^{-1}o)-\rho_{\xi}^f(o,\gamma^{-1}o)$ depend only on the equivalence class $[f]$, but neither on $o$ nor on $f$. 
More precisely, if $o'$ is another point, we have: 
$$
\forall \xi\in\Lambda,\quad \frac{d\nu_o^f}{d\nu_{o'}^f}(\xi)\,=\,\exp(\delta^f\beta_{\xi}(o',o)-\rho^f_{\xi}(o',o))\,.
$$
And if $f=f'+c+X.\phi$, up to a multiplicative constant, we have: 
$$
\forall \xi\in\Lambda,\quad \frac{d\nu_o^f}{d\nu_o^{f'}}(\xi)=\exp(\phi(v_{o,\xi}))\,, 
$$
with $v_{o,\xi}$ the vector based at $o$ pointing to $\xi$.
\end{prop}

\noindent
For simplicity, we will write 
$\beta_{\xi}^f(o,\gamma^{-1}o)=\delta^f \beta_{\xi}(o,\gamma^{-1}o)-\rho_{\xi}^f(o,\gamma^{-1}o).$
It is another cocycle on $\Gamma\times\D$ associated to $f$, the Radon-Nikodym cocycle of $\nu_o^f$.
Notice also that $\beta^f=\rho^{\delta^f-f}$, so that $[\beta^f]_{\cal H}=[c_o^f]_{\cal H}$.

\noindent
Conversely, Ledrappier characterized the H\"older cocycles which are the Radon Nikodym cocycles associated to a quasi-invariant measure:
more precisely, he constructed a bijection between the equivalence classes of H\"older potentials on $\tm$, the 
equivalence classes of quasi-invariant measures on $\Lambda\subset\D$ with a H\"older cocycle, and the 
cohomology classes on $\D$ of normalized H\"older cocycles.\\
%

\noindent
Recall that we assume for the moment that $f$ is symmetric.
The second step is:
 
\begin{prop}The Radon measure $\widetilde{m}^f$ on $\Lambda^2\times \R$ defined below is invariant by $\Gamma$ and by the geodesic flow: 
$$
d\widetilde{m}^f(v)=\exp\left(\beta_{v^+}^f(o,v)+\beta_{v^-}^f(o,v)\right)\;d\nu_o^f(v^-)\;d\nu_o^f(v^+)\;ds.
$$ 
So it induces on the quotient a finite measure on $\tm$ (normalized to be a probability) 
with support in the non-wandering set $\Omega$: the equilibrium state $m^f$ associated to $f$ (see \cite{kaim2}, \cite{ham}).
Moreover, this measure is mixing w.r.t the geodesic flow $\Phi$, \cite{mbab2}. 
Finally, it depends only on the equivalence class $[f]$ of $f$.
\end{prop}

\noindent
The quasi-product structure of $\widetilde{m}^f$ suggests to decompose it into two families of measures. 
First, the family $(\mu^f_{H^+})_{H^+\in{\cal H}}$ defined by:
$$
d\mu^f_{H^+(u)}(v)=\exp(\beta^f_{v^+}(o,v))\;d\nu_o^f(v^+)\,,
$$ 
is a $\Gamma$-invariant family of measures on the leaves of the foliation, 
in the sense that $\gamma_*\mu^f_{H^+}=\mu^f_{\gamma H^+}$.

\noindent
Recall that a quasi-invariant transverse measure for a foliation 
is completely determined by its restriction to a family of transversals passing at all points.
Then, we have: \\

\noindent
{\bf Proposition 1'} {\it The transverse measure to the foliation $\widetilde{\cal W}^{su}$ defined on each transversal $T=\widetilde{W}^s(w)$ by: 
$$\forall v\in T,\quad d\mu_T^f(v)=\exp(\beta_{v^-}^f(o,v))\,d\nu_o^f(v^-)\,ds$$
is invariant by $\Gamma$ in the sense that $\gamma_*\mu^f_{T}=\mu^f_{\gamma T}$. So it induces on $\tm$ 
a transverse measure to ${\cal W}^{su}$ (still denoted by $\{\mu_T^f\}$), which is quasi-invariant by holonomy, with respect to the cocycle $\rho^f$. 
}\\

\noindent
The proof of Proposition 1' is straightforward.\\

\noindent
{\bf Remark:} The above transverse measure $\{\mu_T^f\}$ is constructed using the equilibrium state $m^f$ of $f$.
Thus we can make the analogy between our geometric situation and the symbolic one. In the case of the full 
shift $\Sigma=\{1,...,k\}^{\Z}$ endowed with the shift $\sigma$, the unique equilibrium state $m^f$ coincides
with the unique {\em Gibbs measure} associated to $f$, that is the measure satisfying 
$\zeta_*m^f=\exp(\rho^f)\,m^f$ (see \cite{bowen}, \cite{keller}), 
with $\zeta$ a map from a cylinder to another. 
Since a cylinder can be considered as a transversal to the 
strong stable manifolds of $\sigma$, the map $\zeta$ can be viewed as a holonomy map.
Thus Theorem \ref{1bis} gives a similar result in the case of a manifold: 
the unique quasi-invariant transverse measure with respect to $\rho^f$ is "induced" by the unique 
equilibrium state $m^f$, which can therefore be called a Gibbs measure.\\

\noindent

\noindent
In our geometrical frame, recall that there exists a nice set of leaves of the foliation:
the set ${\cal H}\simeq\D\times\R$ of horospheres. 
Better than working  with the family of measures $\{\mu_T\}$ defined above, we would like to obtain a unique measure on 
 ${\cal H}$. To do this, we need  a holonomy invariant family of measures, or 
equivalently measures $\mu_T$ on transversals $T=\widetilde{W}^s(w)$ which do not depend on $T$, i.e. on $v^+$. 
The only term in the expression of $\mu_T$ which is not 
invariant by holonomy is $\exp(-\rho^f_{v^-}(o,v))$ (when $f$ is not a constant).
So 
we divide the measures $\mu_T$ by the density $\exp(-\rho^f_{v^-}(o,v))$, and we multiply the measures 
$\overline{\mu}^f_{H^+ u}$ by the same density. Then we obtain new families of measures: 
$$
d\overline{\mu}^f_T(v)=ds\,d\nu_o^f(v^-)\,\exp\left(\delta^f \beta_{v^-}(o,v)\right)\;\quad\mbox{ and}$$
 $$ d\overline{\mu}^f_{H^+ u}(v)=d\nu_o^f(v^+)\,\exp\left(\beta^f_{v^+}(o,v)-\rho^f_{v^-}(o,v)\right).
$$

\noindent
These families are now quasi-invariant by $\Gamma$. More precisely, the measures on leaves satisfy:
\begin{eqnarray}\label{mufh}
\forall v\in\Lambda^2\times\R,\;\forall \gamma\in\Gamma, \quad  
\frac{d\gamma^{-1}_*\overline{\mu}^f_{H^+(\gamma v)}}{d\overline{\mu}^f_{H^+ v}}(v) = 
\exp ( -c_o^f(\gamma,v^-) ),
\end{eqnarray}
\noindent
Pushing by the flow leads to the property: 
\begin{eqnarray}\label{flo}
\forall v\in\Lambda^2\times\R,\;\forall t\in\R,\quad 
\frac{d\Phi^{-t}_*\bar{\mu}^f_{H^+(\Phi^t v)} }{d\bar{\mu}^f_{H^+ v}}(v)=\exp(t\delta^f).
\end{eqnarray}

\noindent
The  family $\overline{\mu}^f_T$  is also quasi-invariant by $\Gamma$,
with the opposite cocycle: $ c_o^f(\gamma,v^-)=\rho^f_{v^-} (\gamma^{-1} o, o) $.
But it is now invariant by holonomy, so it induces a measure $\chap{\mu}_o^f$ on the quotient on the space ${\cal H}$ of leaves.
More precisely, we have: 

\begin{lem}The measure $d\chap{\mu}_o^f(\xi,s)=\exp(-\delta^f s)\,ds\,d\nu_o^f(\xi)$ on ${\cal H}$ is supported on $\Lambda\times\R$, and
$\Gamma$-quasi-invariant with respect to the cocycle $ c_o^f$. 
\end{lem} 
%

\noindent
Remark that in the case of a constant potential, we did not change anything, 
these two families of measures are the same as in the beginning, and in particular, 
the family $(\mu^f_T)$ is at the same time $\Gamma$-invariant and invariant under the holonomy pseudogroup.

\noindent
As the topological pressure satisfies $\delta^{f+c}=\delta^f+c$, we can assume that $\delta^f>0$. It
 is simply a normalisation convention which allows us
to assume that, when pushed by the flow $\Phi^t$, with $t>0$, 
the family of measures $\bar{\mu}^f_{H^+}$ is expanded by a factor $\exp(t\delta^f)$.

\noindent
Let us now summarize  the main properties of the above measures in a  table.

\begin{center}
\setlength{\tabcolsep}{0.2cm}
\begin{tabular}{|c|c|c|}

\hline
Measures and their  support &  $\Gamma$-Quasi-Invariance  & Q-I by holonomy  \\
\hline
  $\nu_o^f$ on $\Lambda\subset\D$ & $\exp(\beta^f)$& --- \\
\hline
$\mu^f_{H^+}$ on leaves $H^+$ of ${\cal W}^{su}$ & Invariant & ---\\
\hline 
$\bar{\mu}^f_{H^+}$ on  leaves $H^+$ of $\widetilde{\cal W}^{su}$ & $\exp(-c_o^f)$  & --- \\
\hline
$\mu^f_T$ on transversals to ${\cal W}^{su}$ & Invariant & $\exp(\rho^f)$  \\
\hline
$\bar{\mu}^f_T$  on transversals to $\widetilde{\cal W}^{su}$& $\exp(c_o^f)$  & Invariant \\
\hline 
$\chap{\mu}_o^f$ on $\Lambda\times\R\subset {\cal H}$ & $\exp(c_o^f)$& ---\\
\hline
\end{tabular}
\end{center}


\section{Horospherical means, Proofs of Theorems \ref{1bis} and \ref{1bis} bis}\label{equicontinuite}



The aim of this section is the proof of the unicity results of a quasi-invariant measure (Theorems \ref{1bis} and \ref{1bis} bis). 
We use an intermediate  result of equidistribution of horospherical means pushed by the flow $\Phi^t$ when $t\to\infty$ \cite[thm 4]{mbab}, 
that we recall below, and
some properties of continuity of the measure of balls, 
proved in section \ref{par2}. 
A formula (section \ref{autoadjonction}) allows in section \ref{preuves} to prove theorem \ref{1bis} bis. 
Theorem \ref{1bis} is obtained as a corollary of the preceding.  

\noindent
Finally, we give in section \ref{equidunif} an alternative method to prove Theorems \ref{1bis} bis and \ref{convergence} 
in a particular case.
In fact, under the aditional assumption that the boundary of all horospherical balls 
has measure zero, we are able to prove directly a uniform equidistribution property (Theorem \ref{convergence} bis), and 
to deduce easily Theorem \ref{1bis} bis. 
(In the general case, recall that we will deduce in section \ref{equid} Theorem \ref{convergence}
from Theorem \ref{1bis} bis.)\\


\noindent
Given a potential $f$, the mean $M_{r,u}(\psi)$ of a map $\psi:\tm \to\R$ on $B^+(u,r)$ is defined by:
$$
M_{r,u}(\psi)\,=\,\frac{1}{\bar{\mu}^f_{H^+u}(B^+(u,r))}\,\int_{B^+(u,r)}\widetilde{\psi}\,d\bar{\mu}^f_{H^+u}\,,
$$
 with $\widetilde{\psi}$ the $\Gamma$-invariant lift of $\psi$ on $\ttm$. 
Notice that these means cannot be defined directly on $\tm$, since the measure $\bar{\mu}^f_{H^+u}$ is not $\Gamma$-invariant.
However, they are $\Gamma$-invariant, so we shall consider them as probability measures on $\tm$.\\
For simplicity, we denote by  $M_{r,u}^t(\psi):=\Phi^t_* M_{r,u}(\psi)$ the mean pushed by the flow. 
The geometrical property of the horospherical distances then implies
the simple, but fundamental property that the mean of $\psi\circ\Phi^t$ on a ball $B^+(u,r)$ 
equals the mean of $\psi$ on the bigger ball $B^+(\Phi^t u, re^t)$:
\begin{eqnarray}\label{magic}
\forall u\in\ttm,\,\forall r>0,\,\forall t\in\R,\quad M_{r,u}^t(\psi)=M_{re^t,\Phi^t u}(\psi).
\end{eqnarray}

\noindent
Recall also that, from a general result in \cite{mbab}, we have:
\begin{theo}[Babillot, \cite{mbab}]\label{martine}Assume that $\Gamma$ is cocompact or convex-cocompact,
and the geodesic flow is topologically mixing on its nonwandering set $\Omega$. 
For all $u\in\Lambda^2\times\R$, for all $r>0$, the mean
$M_{r,u}^t$ converges weakly to $m^f$ when $t\to +\infty$.
\end{theo}

\noindent 
This result is a consequence of the mixing of the measure $m^f$, 
which is established when the geodesic flow is topologically mixing. 

\noindent
Notice that, after fixing a positive and continuous  map $\psi:\tm\to\R$, 
we will often consider these means $M_{r,u}^t(\psi)$ as {\em functions of the variable} $u\in\Lambda^2\times\R$.
By $\Gamma$-invariance,
we consider them either as functions on
$\Omega$ or $\Lambda^2\times\R$. 

%

%



\subsection{Averaging on horospherical balls}\label{par2}

In this first section, we shall study the map 
$$ (u,r)\mapsto\int_{B^+(u,r)}\psi\,d\bar{\mu}^f_{H^+}.$$ 
It is not continuous, but we shall prove that it has a regularity property.
The main result is Lemma \ref{116}.
This paragraph is inspired from the sequence of "technical lemmas" of Roblin, \cite{roblin}, paragraph 1H.

\noindent
If $u\in\ttm$ and $v\in H^-(u)$, 
let us introduce the map $P_{u,v}:H^+(u)\to H^+(v)$ (see Figure 4 below),
which sends $w=(u^-,w^+,s(u))\in H^+(u)$
to the intersection $w'$ of the geodesic $(v^-,w^+)$ with $H^+(v)$. In other words, 
$w'=(v^-,w^+,s(v))\in H^+(v)$ (when it makes sense, i.e. when $w^+\neq u^-,v^-$). 
 \begin{figure}[!ht]
\begin{center}
\input{puv.pstex_t}
\caption{The map $P_{u,v}:w\to w'$}
\end{center}
\end{figure}

\noindent
In the following, we prove some lemmas which express mainly three important properties; 
first, these maps $P_{u,v}$ 
are uniformly close to the identity when $v$ is close to $u$ (Lemma \ref{13}).
Second, all these measures $\overline{\mu}^f_{H^+}$ have a continuous (and positive) Radon-Nikodym derivative
with respect to the measure 
$\nu_o^f$ on $\D$, which allows to compare the measures $\overline{\mu}^f_{H^+(v)}$ and ${P_{u,v}}_*\overline{\mu}^f_{H^+(u)}$: 
they are equivalent with continuous 
Radon-Nikodym derivative (Lemma \ref{14}).
Third, the flow contracts the stable manifolds. 
All these properties are proved in some preliminary lemmas, and collected in Lemma \ref{116},
which  will allow to prove Theorem \ref{1bis} bis.

\noindent
To make the discussion more precise, we have to introduce neighborhoods, or "cells", in which we will let $v$ vary: 
if $u\in\ttm$, and $r_1,r_2,r_3$ are three positive numbers, we note 
$$
C(u,r_1,r_2,r_3)=\cup_{|s|<r_3} \Phi^s \cup_{v_1\in B^-(u,r_1)}P_{u,v_1} (B^+(u,r_2)).
$$

\begin{figure}[!ht]
\begin{center}
\input{cellule.pstex_t}
\caption{Cells}
\end{center}
\end{figure}
\noindent 
If $\psi$ is a positive Borel function on $\ttm$, we set:\\
$\displaystyle \psi_{\varepsilon}(w)=\sup\{\psi(v),\,v\in\cup_{|s|<\varepsilon}\Phi^s B^-(w,\varepsilon)\}$,  and\\
$\displaystyle \psi_{-\varepsilon}(w)=\inf\{\psi(v),\,v\in\cup_{|s|<\varepsilon}\Phi^s B^-(w,\varepsilon)\}$.
%
%
Then, if $A\subset\ttm$ is a Borel set, we have $[\un_A]_{\pm\varepsilon}=\un_{A_{\pm\varepsilon}}$, 
and for all $\varepsilon>0$, and $t\geq 0$, 
$[\Phi^{-t}A]_{\varepsilon}\subset\Phi^{-t} (A_{\varepsilon})$, and $ [\Phi^{-t}A]_{-\varepsilon}\supset\Phi^{-t} (A_{-\varepsilon}).$
This means simply that if $t\geq 0$, the flow $\Phi^t$ contracts the stable manifolds, and expands 
the unstable manifolds.\\

\noindent
In the first lemma, we compare balls $B^+(u,r)$ and $B^+(v,r)$ when $v\in H^-(u)$.

\begin{lem}[Roblin, \cite{roblin}]\label{13}  Let $K\subset\ttm$ be a compact, and $\varepsilon>0$. 
There exists $\delta>0$ such that for all $u\in K$, if
$v\in B^-(u,\delta)$ and $w\in B^+(u,3)$, then
$P_{u,v}w \in\Phi^s B^-(w,\varepsilon)$, with $|s|<\varepsilon$, and if  $r\in[1,2]$,
$ B^+(v,re^{-\varepsilon})\subset P_{u,v}(B^+(u,r))\subset B^+(v,re^{\varepsilon}).$
\end{lem}

\noindent
In other words, the map $P_{u,v}:B^+(u,3)\to H^+(v)$ restricted to the ball $B^+(u,3)$ in $H^+(u)$, 
is uniformly closed to the identity of $B^+(u,3)$ (up to $\exp(\pm\varepsilon)$)
for the topology of uniform convergence on compact sets.
\begin{figure}[!ht]
\begin{center}
\input{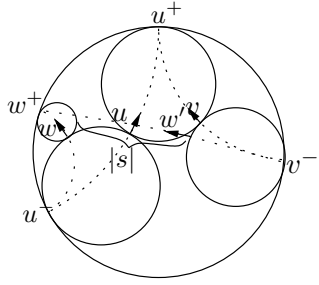}
\caption{Cross ratio of $( u^-,v^- ,u^+,w^+)$ }
\end{center}
\end{figure}

\begin{demo}The assumptions of the Lemma mean that if $\delta$ is small enough, 
the $4$-uple $(u^-, u^+, v^-, w^+)$ varies in a compact set of $\DD\times\DD$. One easily checks that 
the quantity $s$ such that  $P_{u,v}w \in\Phi^s H^-(w)$ is the cross ratio 
$B(v^-,u^-,u^+,w^+)="d(v^-,u^+)+d(u^-,w^+)-d(v^-,w^+)-d(u^-,u^+)"$ of the four points. 
Moreover, the cross ratio is a continuous map, \cite{Otal}, which vanishes if $u^-=v^-$ (or $u^+=w^+$). 
Thus if $\delta$ is small enough, $v^-$ is uniformly closed to 
$u^-$, which gives $P_{u,v}w \in\Phi^s H^-(w)$. Moreover, an easy computation shows that\\
$\displaystyle \frac{d_{H^+}(P_{u,v}w,v)}{d_{H^+}(u,v)}=e^{s/2}$, which concludes the proof. 
\end{demo}

\noindent
The following lemma is a corollary of the above and of continuity of densities of measures 
$\overline{\mu}_{H^+ u}^f$ 
on horospheres with respect to the measure $\nu_o^f$ on the boundary.

\begin{lem}\label{14} Let  $K\subset \tm$ be a compact, and $\varepsilon>0$.
There exists $\delta>0$, such that for all $u\in K$, and $v\in B^-(u,\delta)$, 
for all $w\in B^+(u,3)$ and $w'=P_{u,v}w$, the following quantity is uniformly close to $1$ (up to $\exp(\pm\varepsilon)$):
$$
\frac{d\overline{\mu}_{H^+ v}^f}{d(P_{u,v}*\overline{\mu}_{H^+ u}^f)}(w')=
\exp\left(\beta^f_{w^+}(w,w')-\rho^f_{v^-}(o,w')+\rho^f_{u^-}(o,w)\right).
$$
\end{lem}

\begin{demo}Recall first that 
$\beta^f_{w^+}(w,w')-\rho^f_{v^-}(o,w')+\rho^f_{u^-}(o,w)=\delta^f \beta_{w^+}(w,w')-\rho^f_{w^+}(w,w')-\rho^f_{v^-}(o,w')+\rho^f_{u^-}(o,w)$. 
We have $|\beta_{w^+}(w,w')|\leq\varepsilon$, after \ref{13}. 
We have also 
$|\rho^f_{w^+}(w,w')|\leq C(f) \,D(w,w')^{\alpha(f)}$, 
because $f$ is H\"older, which gives the desired upper bound by \ref{13}. 
Finally, the last term can be written as $\exp(A(P_{u,v}w)-A(w))$, with   
$A(w)=\rho^f_{w^-}(o,w)$, and the result follows from the continuity of $A$ on $\ttm$. Indeed, set $w^t=\Phi^{-t}w$, 
$y\in W^{su}(u)$ the vector of $W^{su}(u)$ 
tangent to the geodesic $(u^-,o]$, and   
$y^t=\Phi^{-t}y$. 
As $f$ is H\"older, and the curvature is bounded above by a negative constant,
there exists a $T>0$, uniform on $K$, such that if $t\geq T$, 
$A(w)$ is $\varepsilon$-closed to $\int_o^{y^t}f-\int_{w}^{w^t}f$ . But for all fixed  $t>0$, 
this quantity is continuous, hence uniformly continuous in $w\in K$. So $A$ is continuous, and the result follows.
\end{demo}

\noindent
We eventually arrive at the last lemma, which summarizes the main properties of the measure $\overline{\mu}^f_{H^+(u)}$ on balls.

\begin{lem}\label{116}Let  $\psi$ be a positive Borel function on $\ttm$, $\varepsilon>0$ (small),  and $K\subset\tm$ a compact set. There exists
 $r_1,r_2,r_3>0$ small enough, so that: for all $u\in K$, for all 
 $v\in C(u,r_1,r_2,r_3)$, $r\in[1,2]$, we have:
$$
e^{-\varepsilon}\int_{B^+(u,re^{-\varepsilon})}\psi_{-\varepsilon}\,d\overline{\mu}^f_{H^+(u)}
\leq \int_{B^+(v,r)}\psi\,d\overline{\mu}^f_{H^+(v)}
\leq e^{\varepsilon} \int_{B^+(u,re^{\varepsilon})}\psi_{\varepsilon}\,d\overline{\mu}^f_{H^+(u)}.
$$
\end{lem}

\begin{dem}The proof follows closely that of Roblin, \cite{roblin}, and goes as follows. 
Let $\varepsilon>0$, small enough so that $r e^{3\varepsilon}\leq 2e^{3\varepsilon}\leq 3$. 
Fix $u\in K$, and $v\in C(u,r_1,r_2,r_3)$. By definition of cells, 
we can write $v=\Phi^s v_2$, with $|s|<r_3$, $v_2\in P_{u,v_1} B^+(u,r_2)\subset H^+ (v_1)$ and $v_1\in B^-(u,r_1)$.
First, we push the ball $B^+(v,r)$ by $\Phi^{-s}$ to the ball $B^+(v_2,r e^{-s})$, with $s$ very small.
By relations (\ref{flo}) and (\ref{magic}), the integrals of $\psi$ on these balls are almost the same.
Second, we move the ball $B^+(v_2,r e^{-s})$ by the map $P_{u,v_1}^{-1}$. Lemma \ref{13} says that the image is uniformly close to the ball 
$B^+(u,re^{-\varepsilon})$ and Lemma \ref{14} says that the derivative of the measure $\mu^f_{H^+}$ is uniformly close to $1$. 
Third, 
the image of $\psi$ by these operations is close to $\psi_{-\varepsilon}$. 
The right inequality is proved in the same way.
\end{dem}

\subsection{An autoadjonction property}\label{autoadjonction}

It is obvious that for a flow $(\Phi^t)_{t\in\R}$ on a compact space $X$, the equality
$$\int_X\frac{1}{T}\int_0^T\psi \circ\Phi^t\,dt d\nu=\int_X\psi\,d\nu$$ 
for all $T>0$ and any $\Phi$-invariant measure $\nu$ allows to prove unique ergodicity if we know that 
all Birkhoff averages converge to a constant.

\noindent
In our situation, the following equality allows to do the same trick, 
replacing Birkhoff averages by horospherical means: we will use it to prove unique ergo-dicity, 
using a weaker property than equidistribution of horospherical means.

%
%
%
\begin{lem}\label{3}Let $\chap{\nu}$ be any fixed quasi-invariant measure on
 ${\cal H}$ with the cocycle $c_o^f$,
and $M$ the measure on $\DDR$ defined by
$dM(v)=d\chap{\nu}(H^+)\,d\mcb(v)$. If $\psi$ a $\Gamma$-invariant positive measurable map on $T^1\M$, 
and $D$ a measurable fundamental domain for the 
action of $\Gamma$ on $T^1\M$, then:
\begin{eqnarray*}
\int_D\,dM(u) \,M_{r,u}(\psi)\,
=\int_D\,dM(v)\psi(v)\int_{B^+(v,r)}\frac{d\overline{\mu}^f_{H^+(v)}(u)}{\overline{\mu}^f_{H^+(v)}(B^+(u,r))}\;.
\end{eqnarray*}
\end{lem}

\noindent
Before we proceed to the proof, let us make a remark: 
if we knew that uniform equidistribution of horospherical means 
when $r\to\infty$ holds, then, as in the case of a flow, 
the above formula would allow to deduce the unique ergodicity result from the equidistribution property.
In fact, we will use this approach in a particular case in section \ref{equidunif}.


\begin{demo}By definition, $dM=d\chap{\nu}\,\,d\mcb$, and
the left term can be rewritten:
$$
\int_{{\cal H}}d\chap{\nu}(H^+(u))\int_{H^+(u)\times H^+(u)}
d\mcb(u)\,d\mcb(v) \frac{{\un}_D(u)\,{\un}_{B^+(u,r)}(v)\,\psi(v)}{\overline{\mu}^f_{H^+(u)}(B^+(u,r))}.
$$
The proof is then based on the fact that $\psi$ is
 $\Gamma$-invariant, $\mcb$ and 
$\chap{\nu}$ are $\Gamma$-quasi-invariant (resp. with cocycles $-c_o^f$ and $+c_o^f$)
 and on two observations: first ${\un}_{B^+(u,r)}(v)$ 
is symmetric in $u$ and $v$; 
second, $B^+(u,r)=\sqcup_{\gamma\in\Gamma} \,B^+(u,r)\cap \gamma D$. \\
If we set $v=\gamma v'$ and $u=\gamma u'$, we have: 
${\un}_{B^+(u,r)}(v)=\sum_{\gamma\in\Gamma}{\un}_{B^+(u',r)}(v') {\un}_{D}(v')$.\\
If we decompose the above integral into a sum over 
 $\gamma\in\Gamma$, we find that it equals: \\

\begin{eqnarray*}
\sum_{\gamma\in\Gamma} \int_{{\cal H}}d\chap{\nu}(\gamma H^+(u'))
\int_{\gamma H^+(u')\times \gamma H^+(u')}d\overline{\mu}_{\gamma H^+(u')}^f(\gamma u')
\,d\overline{\mu}_{\gamma H^+(u')}^f(\gamma v')\,\, \times\\
\times\,\,\frac{{\un}_{\gamma^{-1}D}(u') {\un}_{B^+(u',r)}(v'){\un}_D(v')\psi((\gamma v')}{\overline{\mu}^f_{\gamma H^+(u')}(\gamma B^+(u',r))}
\end{eqnarray*}
We then use the various (quasi-)invariance relations by $\Gamma$;
we remark also that 
$\sum_{\gamma}{\un}_{\gamma^{-1}D}(u')=1$, and $H^+(u)=H^+(v)$, and we get:
$$
\int_{{\cal H}}d\chap{\nu}(H^+(v))\int_{H^+(v)\times H^+(v)}d\mcb(v)\,d\mcb(u)\,\frac{\un_D(v)
{\un}_{B^+(v,r)}(u)\psi(v)}{\overline{\mu}^f_{H^+(v)}(B^+(u,r))}.
$$
This expression is the same than the first one, with the roles of  $u$ and $v$ exchanged.
\end{demo}


\subsection{Proof of Theorems \ref{1bis}  and \ref{1bis} bis}\label{preuves}

%
%
%

We prove first Theorem \ref{1bis}  bis.

\begin{demon}{of Theorem \ref{1bis} bis} Let $[c]_{\cal H}$ be a cohomology class of H\"older cocycles for the action of $\Gamma$ on ${\cal H}$. 
By assumption, this class contains a cocycle on $\D$, so we can assume that $c$ is this cocycle. In \cite{Ledrappier}, Ledrappier showed that every H\"older cocycle
on $\D$ is of the form $c=\lambda \beta^F$, 
with $F$ a H\"older potential and $\lambda\in\R$, 
whence $c=\lambda\rho^{\delta^F-F}=\rho^{\lambda\delta^F-\lambda F}$. 
Finally, $c$ can be written $c=c_o^f$, with $f$ a H\"older 
potential.

\noindent
Let $\chap{\nu}$ 
be a nonzero $\Gamma$-quasi-invariant measure with the cocycle $c_o^f$, 
and $M$ the measure on 
$\DDR$ defined by:
$dM(v)=d\chap{\nu}(H^+(v))\,d\mcb(v)$. Let $\psi$ be a positive continuous map on $\tm$, and $\widetilde{\psi}$ 
its lift on $\ttm$. 
Lemma \ref{3} gives:
$$
\int_D\,dM(u) M_{r,u}^0(\psi)
=\int_D\,dM(v)\widetilde{\psi}(v)\int_{B^+(v,r)}\frac{d\overline{\mu}^f_{H^+(v)}(u)}{\overline{\mu}^f_{H^+(v)}(B^+(u,r))}\;.
$$
%
\noindent
Since $\Omega$ is a compact set, we can choose a relatively compact fundamental domain $D$ for the action
of $\Gamma$ on $\Lambda^2\times\R$, with $M(\partial D)=\widetilde{m}^f(\partial D)=0$.\\
We will prove that for all $r>0$, the right term is smaller than ${\rm cst}.\int\psi\,dM$ (inequality (\ref{majoration})), 
and that there exists $r>0$, such that 
the left term is greater than ${\rm cst}'.\int\psi\,dm^f$, so that $m^f\ll M$ (inequality (\ref{minoration})).
The constants do not depend on $\psi$, so that these two inequalities imply $m^f\ll M$.

\noindent
First, since $\Gamma$ is cocompact or convex-cocompact, the following Vitali property is satisfied
(see \cite{roblin}, Proposition 6.3):  
there exists a uniform integer  $N$ on $C=\Lambda^2\times\R$, such that for all $u\in C$, and all $r>0$,  
\begin{eqnarray}\label{Vitali}
B^+(u,r)\cap \overline{C}\subset \cup_{i=1}^N B^+(u_i,\frac{r}{2}).
\end{eqnarray}
The integral on the ball $B^+(u,r)$ is then smaller than the sum of $N$ integrals on balls $B^+(u_i,\frac{r}{2})$. 
Moreover, triangular inequality gives $B^+(v,r)\cap\bar{C}\supset B^+(u_i,\frac{r}{2})\cap \overline{C}$
if $v\in B^+(u_i,\frac{r}{2})$. 
Hence $\overline{\mu}^f_{H^+(u)}(B^+(v,r))\geq \overline{\mu}^f_{H^+(u)}(B^+(u_i,\frac{r}{2}))$, and finally, 
$$
\forall r>0,\,\forall u\in\Lambda^2\times\R,\quad
\int_{B^+(u,r)}\frac{d\overline{\mu}^f_{H^+(u)}(v)}{\overline{\mu}^f_{H^+(u)}(B^+(v,r))}\leq N.
$$
We deduce that:
\begin{eqnarray}\label{majoration}
\forall r>0,\quad 
\int_D\,dM(v)\widetilde{\psi}(v)\int_{B^+(v,r)}\frac{d\overline{\mu}^f_{H^+(v)}(u)}{\overline{\mu}^f_{H^+(v)}(B^+(u,r))}
\leq N\,\int_{\tm}\psi\,dM.
\end{eqnarray}

\noindent
We shall now prove the following inequality:
\begin{eqnarray}\label{minoration}
\exists r>0, \exists C>0, \int_D\,dM(u) M_{r,u}^0(\psi) \geq C\,\int_{\tm}\psi\,dm^f.
\end{eqnarray}

\noindent
Let $0<\varepsilon<\log 2$ be small enough so that $\int_{\tm}\psi_{-\varepsilon}dm^f\ge \frac{1}{2}\int_{\tm}\psi\,dm^f$.
We can cover $D$ by a finite number of cells $C(u_i)$, $1\le i\le k$, given by Lemma \ref{116} for the above $\varepsilon$. 
If $v\in C(u_i)$, this Lemma gives for all $t\ge 0$, 
$$
M_{1,v}^t(\psi)\ge e^{-2\varepsilon}\,
M_{e^{-\varepsilon},u_i}^t(\psi_{-\varepsilon})\,\frac{\bar{\mu}^f_{H^+}(B^+(u_i,e^{-\varepsilon}))}{\bar{\mu}^f_{H^+}(B^+(u_i,e^{\varepsilon}))}
$$
Moreover, $M_{e^{-\varepsilon},u_i}^t(\psi_{-\varepsilon})\to\int_{\tm}\psi_{-\varepsilon}dm^f$
when $t\to +\infty$ (Theorem \ref{martine}), whence:  
$$\exists T>0,\: \forall u_i,\: M_{e^{-\varepsilon},u_i}^T(\psi_{-\varepsilon})\ge \frac{1}{2} \int_{\tm}\psi_{-\varepsilon}dm^f.$$ 
We deduce that for all $v\in D$ (with $\varepsilon<\log 2$):
$$
M_{1,v}^T(\psi)\ge\left(\frac{1}{16}\int_{\tm}\psi\,dm^f\right)\,\inf_{w\in D}\frac{\bar{\mu}^f_{H^+}(B^+(w,\frac{1}{2}))}{\bar{\mu}^f_{H^+}(B^+(w,2))}
$$
By $\Gamma$-invariance of the means, the above inequality is true for all $v\in\Lambda^2\times\R$.

\noindent
Now, if $u\in\Lambda^2\times\R$, we know that $M_{e^T,u}(\psi)=M_{1,\Phi^T u}^T(\psi)$ and above inequality applies with $v=\Phi^T u$.
At last, we get $$M_{e^T,u}^0(\psi)\,dM(u)\,\ge c\,\int_{\tm}\psi dm^f,\quad\mbox{with }\quad
c=\frac{1}{16}\,\inf_{w\in D}\frac{\bar{\mu}^f_{H^+}(B^+(w,\frac{1}{2}))}{\bar{\mu}^f_{H^+}(B^+(w,2))}>0.$$ 
Finally, integrating with respect to $M$ on $D$ leads to  (\ref{minoration}) with $C=M(D)c$.\\
This constant $C$ does not depend on $\psi$, 
whence $m^f\ll M$.
We deduce from the above that $\chap{\mu}_o^f\ll \chap{\nu}$.
Indeed, let $A\subset{\cal H}$ be a Borel set, and assume that $\chap{\nu}(A)=0$. Let $0<c<C$ be two constants.
If $A$ is small enough, we can choose for all $H\in A$ a distinguished vector
$v\in H$, and a ray $r$ chosen such that for all $H\in A$, $0<c\leq\overline{\mu}^f_{H^+}(B^+(v,r)\leq C$. 
Then 
\begin{eqnarray*}
\int_A d\chap{\nu}(H)\int_H {\bf 1}_{B^+(v,r)}\leq C\chap{\nu}(A)=0, \quad\quad \mbox{whence}\\
 0\leq c\chap{\mu}_o^f(A)\leq \int_A d\chap{\mu}_o^f(H)\int_H {\bf 1}_{B^+(v,r)}\,d\overline{\mu}^f_{H^+ v}=0.
\end{eqnarray*} 

\noindent
If $A$ is not small enough, we can write it as a countable union of smaller sets.
So we proved $\chap{\mu}_o^f\ll \chap{\nu}$.

\noindent
Let us now prove that $\chap{\mu}_o^f$ is ergodic: 
let $A\subset {\cal H}$ be a Borel set such that for all  $\gamma\in\Gamma$, 
$\chap{\mu}_o^f(\gamma A\triangle A)=0$, and $\chap{\mu}_o^f(A)>0$. 
Then the measure $\chap{\mu}^f_{o|A}$ is nonzero and quasi-invariant with cocycle $-c_o^f$.
By what precedes, it satisfies  $\chap{\mu}_o^f\ll\chap{\mu}^f_{o|A}$, and this implies $\chap{\mu}_o^f(A^c)=0$. 
So we have ergodicity.

\noindent
Now let $\chap{\nu}$ be a nonzero and quasi-invariant measure with cocycle $c_o^f$
and let us show that it is also ergodic: 
if  $A$ is a Borel set in  ${\cal H}$, $\Gamma$-invariant $\chap{\nu}$-a.e, and if $\chap{\nu}(A)>0$, we can decompose the measure  
$\chap{\nu}$ into $\chap{\nu}=\chap{\nu}_{|A}+\chap{\nu}_{|A^c}$. By the above, there exist measurable maps $\phi$ and 
$\phi'$ such that $\chap{\mu}_o^f=\phi\,\chap{\nu}=\phi\, \chap{\nu}_{|A}+\phi\,\chap{\nu}_{|A^c} $ and  $\chap{\mu}_o^f=\phi'\,\chap{\nu}_{|A}$.
We deduce that $\chap{\mu}_o^f(A^c)=0$, and then necessarily that  $\chap{\nu}_{|A^c}=0$, which means that $\chap{\nu} $ is ergodic.
 
\noindent
Moreover, the Radon-Nikodym derivative $\frac{d\chap{\mu}_o^f}{d\chap{\nu}}$ is  $\Gamma$-invariant
 $\chap{\nu}$-a.e, hence constant by ergodicity of $\chap{\nu}$. So there exists a unique  
(up to a multiplicative constant) measure which is
 $\Gamma$-quasi-invariant with cocycle $c_o^f$.

\noindent
Recall now that we consider the cohomology class on ${\cal H}$ of the H\"older cocycle $c$ on $\D$.
Let  $c'\in[c]_{\cal H}=[ c_o^f]_{\cal H}$. There exists a positive Borel map 
$\varphi$ on ${\cal H}$, such that for all $\gamma\in\Gamma$ and $(\xi,s)\in{\cal H}$,
$c'(\gamma,(\xi,s))=c_o^f(\gamma,\xi)+\varphi(\gamma(\xi,s))-\varphi(\xi,s)$.
If $\chap{\nu}'$ is a Borel Radon measure on 
${\cal H}$, which is $\Gamma$-quasi-invariant with the cocycle $c'$
then an easy calculation shows that $\exp(-\varphi)\chap{\nu}'$ 
is quasi-invariant  with the cocycle $c_o^f$, it is then unique (up to a constant), 
and equals $\chap{\mu}_o^f$,
so $\chap{\nu}'$ is also unique, and equals $\exp(\varphi)\chap{\mu}_o^f$.
\end{demon}

\begin{demon}{of Theorem \ref{1bis}}It is clear that there exists a bijection between 
 measures transverse to the strong unstable foliation on $T^1 M$, which are
quasi-invariant with the cocycle $\rho^f$, and measures transverse to the strong unstable foliation on
$T^1 \M$, which are $\Gamma$-invariant and quasi-invariant by holonomy with cocycle $\rho^f$.

\noindent
Let us consider the map: $M_T\to N_T$, which sends a $\Gamma$-invariant transverse measure, quasi-invariant by holonomy with cocycle $\rho^f$, 
to the measure $dN_T(v)=\exp(\rho^f_{v^-}(o,v))\,dM_T(v)$. It is clearly a bijection.
An easy computation shows that $N_T$ is a measure transverse to $\widetilde{\cal F}$
which is holonomy invariant and $\Gamma$-quasi-invariant with cocycle $c_o^f$. 
We can associate to $N_T$ 
 a $\Gamma$-quasi invariant measure on ${\cal H}$  (with cocycle $c_o^f$), and this correspondance is 1-1.
The existence of a unique quasi-invariant measure on ${\cal H}$ with cocycle $c_o^f$ is then equivalent to the existence of 
a unique measure  transverse to ${\cal F}$ with the cocycle $\rho^f$.

\noindent
Theorem \ref{1bis} then follows from Theorem \ref{1bis}  bis.
\end{demon}


\subsection{Equidistribution of means in a particular case}\label{equidunif}

Lemma \ref{116} means that the map $(u,r)\to \int_{B^+(u,r)}\psi\,d\bar{\mu}^f_{H^+}$ is "almost continuous". 
In this section, we shall see that under an additional assumption, it is really continuous, 
so that we can prove Theorem \ref{convergence} directly, and deduce Theorem \ref{1bis} bis as a corollary.
The assumption is the following:
$$\forall u\in\Lambda^2\times\R,\quad\bar{\mu}^f_{H^+}(\partial B^+(u,1))=0.$$
Note that it is true in the constant curvature case, and in the case of a surface, since the Patterson measure has no atoms.

\noindent
With this assumption, 
we can prove that the means $(M^t_{1,u}(\psi))_{t\ge 0}$ are equicontinuous in $u\in\Lambda^2\times\R$.
We deduce that the convergence in Theorem \ref{martine} is uniform in 
$u\in\Lambda^2\times\R$, and at last, 
relation (\ref{magic}) leads directly to the uniform convergence of the means $M_r(u)$ to $\int_{\tm}\psi\,dm^f$ 
when $r\to\infty$.\\
Note that these two arguments of equidistribution of the means when $t\to +\infty$,
and equicontinuity in $t\ge 0$ were already used in \cite{EP} in an algebraic framework.

\noindent
Under the above assumption, Lemma \ref{116} leads to:

\begin{lem}\label{equi}If for all $ u\in\Lambda^2\times\R,\,\bar{\mu}^f_{H^+}(\partial B^+(u,1))=0$, then 
for all  $u\in\Lambda^2\times\R$, the family  $(M_{1,u}^t(\psi))_{t\geq 0}$ is equicontinuous in $t\geq 0$.
\end{lem}

\begin{dem}Let $\varepsilon_0>0$, we shall find for all $u\in\Lambda^2\times\R$
a neighbourhood $V=V(u,\varepsilon_0)$ of $u$, such that for all $v\in V$, for all $t\geq 0$,
$\displaystyle |M_{1,v}^t(\psi)-M_{1,u}^t(\psi)|\leq \varepsilon_0$.\\
Since $\overline{\mu}^f_{H^+}(\partial B^+(u,1))=0$, 
there exists  a $0<\varepsilon\leq \varepsilon_0$, 
such that :
\begin{eqnarray}\label{bidule}
\left|\frac{\overline{\mu}^f_{H^+}(B^+(u,e^{\varepsilon}))}{\overline{\mu}^f_{H^+}(B^+(u,e^{-\varepsilon}))}\right|
\leq e^{\varepsilon_0}.
\end{eqnarray}
\noindent
Since $\psi$ is continuous on $\Omega$, one may choose 
$\varepsilon$ small enough so that $\psi_{\varepsilon}\leq \psi+\varepsilon_0$.
We put then $V(u,\varepsilon)=C(u,r_1,r_2,r_3)$: the cell given by Lemma \ref{116} for $\varepsilon$.
Now, relation (\ref{bidule}) and Lemma \ref{116}, together with  the fact that 
$[\psi\circ\Phi^t]_{\varepsilon}\leq\psi_{\varepsilon}\circ\Phi^t$ (since the geodesic flow expands unstable manifolds) 
allow to prove that for all $v\in V(u,\epsilon)$, and for all $t>0$, $|M_{1,v}^t(\psi)-M_{1,u}^t(\psi)|\le \varepsilon_0$.
\end{dem}

\noindent
With this Lemma, we can now prove the following equidistribution result for the horospherical means: \\

\noindent
{\bf Theorem \ref{convergence} bis}: {\it Assume that for all $u\in\tm$ $\displaystyle \bar{\mu}^f_{H^+}(\partial B^+(u,1))=0$, then for all 
continuous map $\psi:\tm\to\R$, the 
means $M_{r,u}(\psi)$ converge uniformly in $u\in\Lambda^2\times\R$ to $\int_{\tm}\psi\,dm^f$.}

\begin{demo}We consider now the means as functions on $\tm$. 
By the above Lemma, for all $u\in\Lambda^2\times\R$, 
the means $(M_{1,u}^t)_{t\geq 0}$ are equicontinuous, whence (Theorem \ref{martine}) it is easy to prove that they converge uniformly on
$\Omega$ to $\int_{\tm}\psi\,dm^f$ when $t\to +\infty$.
Let $\varepsilon>0$, there exists a $T>0$, such that:
$$\forall\,  t\geq T,\;\forall \, u\in\Omega, \quad
\left| M_{r,u}^t(\psi)-\int_{\tm}\psi\,dm^f\right|\leq\varepsilon.
$$
\noindent
Let now $v\in\Omega$ be a vector, and $t\geq T$; as $\Omega$ is invariant by the flow, 
the above inequality applies for $u=\Phi^{-t} v$ 
and above relation (\ref{magic}) gives the desired result:
$$
\left| M_{r e^t,v}^0(\psi)-\int_{\tm}\psi\,dm^f\right|\leq\varepsilon.
$$ 
\end{demo}

\noindent
We can now easily deduce Theorem \ref{1bis} bis from this result:
\begin{dem}We use the equality of Lemma \ref{3}. 
By the uniform equidistribution property of means, 
the left term converges to $M(D)\int_{\tm}\psi\,dm^f$, 
and the right term is smaller than $N\,\int_{\tm}\psi\,dM$. 
Then we have directly $m^f\ll M$, and we conclude the proof as in the general case.  
\end{dem}


\section{Growth of leaves of the horospherical foliation, equidistribution of horospheres}\label{equid}

The aim of this section is to prove Theorem \ref{convergence}. 
In fact, using Theorem \ref{1bis} bis, we shall prove a more general result of equidistribution of sets on leaves satisfying a 
certain condition of growth, denoted by $(*)$ (Theorem \ref{p7}).
We deduce Theorem \ref{convergence} as a 
corollary, 
since the horospherical balls have a polynomial growth (Lemma \ref{croissance}), whence they satisfy $(*)$ (Lemma \ref{*}).\\

\noindent
This paragraph is inspired by Plante \cite{Plante}, and Goodman-Plante \cite{Plante2}.
These two references define a notion of ``averaging sets'' for the foliation, i.e. sequences of sets on transversals 
to the foliation, defined  so that they become equidistribued to a holonomy invariant measure. 
The idea here is to introduce an analogous notion to obtain  equidistribution to the holonomy quasi-invariant measure $\{\mu^f_T\}$.
But it is more appropriate to work with the holonomy invariant measure $\{\bar{\mu}^f_T\}$ on transversals to $\widetilde{\cal W}^{su}$. 
The definition of an averaging set in \cite{Plante2} involves the holonomy pseudogroup of the foliation. In our case, this 
pseudogroup is trivial on $\widetilde{\cal W}^{su}$. Moreover, the above transverse measure is $\Gamma$-quasi-invariant, so
we will introduce a condition dealing with $\Gamma$.\\
Let us now fix a sequence $(E_n)_{n\in\N}$ of compact sets of a horosphere $H^+$.
If they satisfy the condition $(*)$ described below, then we shall prove (Theorem \ref{p7}) that
the means on the $E_n$ w.r.t. the measure $\bar{\mu}^f_{H^+}$ 
become equidistribued to $m^f$. 

\noindent
Let $D$ be any proper fundamental domain for the action of $\Gamma$ on $\ttm$, that is an open set $D$,
such that $D=\stackrel{\circ}{\bar{D}}$,
$\cup_{\gamma\in\Gamma}\gamma \bar{D}=\ttm$, and the $(\gamma D)_{\gamma\in\Gamma}$ are pairwise disjoint.
For fixed $n\in \N$, we distinguish three types of elements $\gamma\in \Gamma$:\\

\noindent
{\bf 1- }First, the set $\Gamma_n^1(D)$ of elements $\gamma$ such that $E_n$ ``crosses completely'' $\gamma D$, i.e.
 $$E_n\cap \gamma D=H^+(u)\cap \gamma D;$$
{\bf 2- }Second,  the set $\Gamma_n^2(D)$ of elements $\gamma$ such that $E_n$ crosses partly $\gamma D$, or again
$$\emptyset\neq E_n\cap \gamma D \varsubsetneq H^+(u)\cap \gamma D$$
{\bf 3-} Third, the set $\Gamma_n^3(D)$ of elements $\gamma$ such that $\displaystyle \gamma D\cap E_n=\emptyset$. \\

\noindent
Compacity of the $E_n$ implies finiteness 
of the sets $\Gamma_n^1(D)$ and $\Gamma_n^2(D)$.\\
Condition $(*)$ is then the following: there exists a fundamental domain $D$, such that, with the above notations: 
\begin{eqnarray*}
\quad\quad\quad\quad\quad\quad\quad\quad\lim_{n\to\infty}\,\frac{1}{\bar{\mu}^f_{H^+}(E_n)}\sum_{\gamma\in\Gamma^2_n(D)}\bar{\mu}^f_{H^+}(E_n\cap\gamma D)
\,=\,0\;\quad\quad\quad\quad\quad\quad\quad\;(*)
\end{eqnarray*}

\noindent
In other words, we can forget the  "boundary terms" $\Gamma_n^2.D$:
only the images $\gamma D$ of $D$ completely crossed by $E_n$ must be taken into account.\\
We can now state the following:

\begin{theo}\label{p7}Let $(E_n)_{n\in\N}$ be a sequence of compact sets on the horosphere  $H^+(u)$, with
 $u\in\Lambda^2\times\R$ and $\cup E_n=H^+(u)$. 
If the sequence $(E_n)_{n\in\N}$ satisfy $(*)$, then the sequence  
$$
dM_n(w)\,=\,\frac{1}{\bar{\mu}^f_{H^+}(E_n)}\;{\bf 1}_{E_n}(w)\;d\bar{\mu}^f_{H^+}(w)
$$
converges weakly to $m^f$ when $n\to +\infty$.
\end{theo}

\noindent
We will prove this result at the end of this section. 
Let us first apply it to sequences of balls $B^+(u,r)$, when $r\to\infty$. 
We shall first prove that balls satisfy $(*)$, 
and thus get the equidistribution result stated in the introduction (Theorem \ref{convergence})
as an immediate corollary of Theorem \ref{p7}. 

\begin{lem}\label{*}The sequence of balls $B^+(u,r)$ satisfies $(*)$.
\end{lem}

\noindent
As a straightforward corollary, Theorem \ref{convergence} is proved by the two above results.
\begin{demo}Let $D$ be any fundamental domain. By compacity of $\bar{D}$,  we can find a constant
 $r_0>0$, such that for all $v\in D$, 
$B^+(v,r_0)\cap D=H^+(v)\cap D$. Now if $\emptyset\neq B^+(u,r)\cap \gamma D\varsubsetneq H^+(u)\cap \gamma D$, 
then triangular inequality gives $B^+(u,r+r_0)\cap \gamma D=H^+(u)\cap \gamma D$. 
Moreover, the same argument gives also  $B^+(u,r-r_0)\cap\gamma D= \emptyset$.
This proves that $B^+(u,r)\cap \Gamma^2_r D\subset B^+(u,r+r_0)\setminus B^+(u,r-r_0)$.
By Lemma \ref{croissance}, it is easy to prove that 
$$
\lim_{r\to +\infty}\,\frac{1}{\bar{\mu}^f_{H^+}(B^+(u,r))}\,\sum_{\gamma\in\Gamma^2_r}\bar{\mu}^f_{H^+}(B^+(u,r)\cap \gamma D)\,=\,0.
$$
\end{demo}

\begin{lem}\label{croissance}For all fixed $u\in\ttm$, $\bar{\mu}^f_{H^+(u)}(B^+(u,r))=O(r^{2\delta^f+4\|f\|_{\infty}})$.
\end{lem}

\begin{demo}If we consider $\Phi^t(u)$ instead of $u$, we can assume that the origin $o$ is
in the interior of the horoball defined by $u^-$ and $u$, or again $\beta_{u^-}(o,u)\leq 0$.
As $\nu_o^f(\D)=1$, and by definition of the measure $\bar{\mu}^f_{H^+(u)}$, we have  
$$
\bar{\mu}^f_{H^+(u)}(B^+(u,r))\leq 
\sup_{v\in B^+(u,r)}\left(\exp(\delta^f\beta_{v^+}(o,v)-\rho^f_{v^+}(o,v)-\rho^f_{v^-}(o,v))\right)\,.
$$ 
Let us first find an upper bound for the quantity $|\rho^f_{v^+}(o,v)+\rho^f_{v^-}(o,v)|$. Since the
curvature of $\M$ is negative,
we can find  a triangle $(a,b,c)$ such that
$a\in [o,v^+)$, $b\in [o,v^-)$, $c\in (v^-,v^+)$, which satisfies: $d(o,a)=d(o,b)$, $\beta_{v^+}(a,c)=0$ and $\beta_{v^-}(b,c)=0$, and also
$d(a,b)\leq \delta$, $d(a,c)\leq \delta$,  $d(b,c)\leq \delta$, with $\delta$ a constant which depends only on the 
upper bound of the sectional curvature of $M$. 
Now, as $f$ is H\"older, with exponent $\alpha(f)$, we can easily check that:
$$|\rho^f_{v^+}(o,v)+\rho^f_{v^-}(o,v)| \leq  2\|f\|_{\infty}\,d(o,a)+C(f,\delta).$$
Moreover, $d(o,a)\simeq d(o,c)\simeq d(o, (v^-,v])\leq d(o,v)$, 
as $o$ is in the interior of the horoball defined by 
$v$ and $v^-$. For the same reason, and by convexity of horoballs, the angle at $v$ between $[o,v]$ and $[v,v^+)$ is greater than $\pi/2$,
whence  $d(o,v)\simeq \beta_{v^+}(o,v)$. 
We deduce 
$$
\bar{\mu}^f_{H^+(u)}(B^+(u,r))\leq Cst. \sup_{B^+(u,r)}\exp((\delta^f+2\|f\|_{\infty}).\beta_{v^+}(o,v)).
$$
As the Busemann cocycle is continuous, 
there exists $\alpha>0$, such that $d_o(u^+,v^+)\leq \alpha$ implies
$|\beta_{v^+}(o,v)|\leq \beta_{u^+}(o,u)+1$. 
Moreover, by definition of $d_{H^+}$,
if $d_o(u^+,v^+)\geq \alpha$, we have the following bound: 
$$
\exp(\frac{1}{2}\beta_{v^+}(o,v))\leq \frac{1}{\alpha}\,d_{H^+}(u,v)\,\exp(-\beta_{u^+}(o,u)/2)\leq C(u).r,
$$
and at last, if $v\in B^+(u,r)$, we have: 
$$e^{(\delta^f+2\|f\|_{\infty}) \beta_{v^+}(o,v)}\leq
\max\left(\,e^{(\delta^f+2\|f\|_{\infty}).(\beta_{u^+}(o,u)+1)},\;\,C(u).r^{2\delta^f+4\|f\|_{\infty}}\right).
$$ 
All these inequalities give the desired comparison (if $u$ is fixed).
\end{demo}

\noindent
Let us now prove Theorem \ref{p7}: 
\begin{demon}{of Theorem \ref{p7}}It is enough to show convergence of means for small Borel sets 
$B\subset D\subset\tm$ which can be written on the form  
$B=T\times b$ in a chart of the foliation, with $T$ a transversal to 
${\cal W}^{su}$, and $b$ a ``piece'' of leaf. 
Lift such a Borel set $B$ into $\widetilde{B}_0=\widetilde{T}_0\times\widetilde{b}_0\subset D$, 
and set $\widetilde{B}=\cup_{\gamma\in\Gamma}\gamma \widetilde{B}_0$. 
We can write again the means on the following form:
\begin{eqnarray*}
 M_n(B)  &=&  \frac{1}
{\bar{\mu}^f_{H^+}(E_n)} \sum_{\gamma\in\Gamma}\int_{H^+(u)} {\bf 1}_{E_n\cap\gamma \widetilde{B}_0}(v)\,d\bar{\mu}^f_{H^+}(v) \\
  & = &\int_{\widetilde{T}_0}\frac{1}{\bar{\mu}^f_{H^+}(E_n)} \sum_{\gamma\in\Gamma} e^{\rho^f_{v^-}(o, \gamma^{-1} o)}
{\bf 1}_{\widetilde{T}_0\cap \gamma^{-1} E_n\neq \emptyset}(v)
\bar{\mu}^f_{H^+}(H^+(v)\cap \widetilde{B}_0)\\
&+ & R(n,\widetilde{T}_0) , \\
  & = & \int_{\widetilde{T}_0}d\nu_{\widetilde{T}_0}^n(v)
\int_{H^+(v)}{\bf 1}_{\widetilde{B}_0}(w)\,d\bar{\mu}^f_{H^+}(w)   + 
 \, R(n,\widetilde{T}_0)\,\, \quad \quad\mbox{   with }\\
d\nu_{\widetilde{T}_0}^n(w) &=&\frac{1}{\bar{\mu}^f_{H^+}(E_n)}\sum_{\gamma\in\Gamma} \exp(\rho^f_{v^-}(o, \gamma^{-1} o))
{\bf 1}_{\widetilde{T}_0\cap \gamma^{-1} E_n\neq \emptyset}(w),
\end{eqnarray*}

\noindent
and the rest $R(n,\widetilde{T}_0)$ corresponds to the $\gamma \in \Gamma$ such that 
$E_n$ crosses only partly the box $\gamma \widetilde{B}_0$; 
so it is bounded by: 
$$
R(n,\widetilde{T}_0)\leq \frac{1}{\bar{\mu}^f_{H^+}(E_n)}\sum_{\gamma\in \Gamma^2_n}\bar{\mu}^f_{H^+}(E_n\cap\gamma D).
$$

\noindent
The sequence $(M_n)$ of probabilities is defined on the compact set $\Omega$, 
so it  has limit points for the weak topology, and after the above computation, 
for any fixed transversal $T$, 
these limit points are in correspondance $1-1$ with those of the sequence 
$(\nu_T^n)$. 

\noindent
Let $M$ be a limit point of measures $M_n$; 
if we consider a subsequence, 
we can assume that the  $M_n$ converge to $M$, 
and the measures $\nu_T^n$ also converge to a measure $\nu_T$. 
This family of measures is invariant by holonomy; 
indeed, if we consider two small transversals $T$ and $T'$ to the foliation, such that
$B=T\times b=T'\times b$ in a chart of the foliation, then 
\begin{eqnarray*}
M(B)\,&=&\,\lim_{n\to\infty} \int_T d\nu_T^n(v)\int_{H^+(v)}{\bf 1}_{\widetilde{B}_0}(w)d\bar{\mu}^f_{H^+}(w)\\
&=&\,\lim_{n\to\infty} \int_{T'} d\nu_{T'}^n(v)\int_{H^+(v)}{\bf 1}_{\widetilde{B}_0}(w)d\bar{\mu}^f_{H^+}(w),
\end{eqnarray*}
which gives invariance by holonomy of the limit measures. 
This family of measures is $\Gamma$-quasi-invariant; indeed, if $g\in\Gamma$, it is easy to compute that: 
\begin{eqnarray*}
d\nu_{gT}^n(gv) &= &
\frac{1}{\bar{\mu}^f_{H^+}(E_n)} \sum_{\gamma\in\Gamma}\exp(\rho^f_{g v^-}(o, \gamma^{-1} o))\,{\bf 1}_{gT\cap \gamma^{-1} E_n}(g v)\\
& = & \exp(\rho^f_{v^-}(g^{-1}o,o))\,d\nu_T^n(v)\,=\,\exp(c_o^f(g,v^-))\,d\nu_T^n(v) .
\end{eqnarray*}

\noindent
At last, we proved that a limit point $M$ of $M_n$ 
induces a transverse measure invariant by holonomy and
quasi-invariant by $\Gamma$ with cocycle $c_o^f$, whence (Theorem \ref{1bis}) it equals the transverse measure 
$(\bar{\mu}_T^f)_T$ induced by $m^f$; 
this shows that the sequence $M_n$ becomes equidistributed to $m^f$.
\end{demon}

\section{Application to coverings}\label{nonsym}


\subsection{Nonsymmetric potential}\label{chech}

If $f$ is not symmetric, we set $\check{f}(v)=f(-v)$, 
and we look again at the product measure
$\nu_o^f\times\nu_o^{\check{f}}$ on $\DD$. 
We are still in the case when $\delta^f<+\infty$, 
and we have $\delta^f=\delta^{\check{f}}$ (see \cite{coudene}).
The equilibrium state $m^f$ is still constructed by means of this product, 
but the density is slightly different \cite{coudene}:
$$
d\widetilde{m}^f(v^-,v^+,s)=
ds\,d\nu_o^f(v^+)\,d\nu_o^{\check{f}}(v^-)\,
\exp\left(\beta^f_{v^+}(o,v)+\beta^{\check{f}}_{v^-}(o,v)\right).
$$
The measures of paragraph \ref{measures} can then be written:  
$$
d\chap{\mu}_o^f(\xi,s)=ds\,d\nu_o^{\check{f}}(v^-)\,\exp(-\delta^f s),\quad\quad\quad\mbox{ and }
$$
$$
d\overline{\mu}^f_{H^+ u}(v)=d\nu_o^f(v^+)\,\exp(\beta^f_{v^+}(o,v)-\rho^{\check{f}}_{v^-}(o,v)).
$$
\noindent
And the relations of $\Gamma$-quasi-invariance become:
$$
\frac{d(\gamma^{-1}_*\chap{\mu}_o^f)}{d\chap{\mu}_o^f}(\xi,s)=\exp(c_o^{\check{f}}(\gamma,\xi)), \quad\mbox{and}\quad
\frac{d(\gamma^{-1}_*\overline{\mu}_{H^+ (\gamma v)}^f)}{d\chap{\mu}_{H^+ v}^f}(v)=\exp(-\rho^{\check{f}}_{v^-}(\gamma^{-1} o,o))
$$
All proofs apply in the same way, 
since they are based on the H\"older continuity of $f$ and geometry of $M$, whence Theorems \ref{1bis} and \ref{1bis} bis are still true.


\subsection{Proof of Theorem \ref{cover}}

The above paragraph applies in particular to the case of a closed $1$-form $\alpha$ on $M$. Indeed, a $1$-form on $M$ 
is a smooth map from $M$ into $T^*M$ which to each $x\in M$ associates a linear form $\alpha_x$ on $T_x M$. 
In particular, it is a map from $TM$ into $\R$, that we can restrict to $\tm$. For all $x\in M$, 
linearity of $ \alpha_x$ gives antisymmetry of $\alpha$: $\alpha(-v)=-\alpha(v)$. 
Moreover, as $\alpha$ is closed, we have 
$\rho_{\xi}^{\check{\alpha}}(x,y)=\int_x^y \check{\alpha}=-\rho_{\xi}^{\alpha}(x,y)$, 
whence (cf \cite{mbab3}):
$$
d\widetilde{m}^{\alpha}(v^-,v^+,s)=
ds\, d\nu_o^{\alpha}(v^+)\, d\nu_o^{\check{\alpha}}(v^-)\,\exp\left(\delta^{\alpha}\beta_{v^-}(o,v)+\delta^{\alpha}\beta_{v^+}(o,v)\right).
$$
The measure $\chap{\mu}_o^{\alpha}$ can be written: 
$d\chap{\mu}_o^{\alpha}(u^-,s)=ds\, d\nu_o^{\check{\alpha}}(v^-)\,\exp(-\delta^{\alpha}s).$
Moreover, since $\alpha$ is closed, 
we have $c_o^{\check{\alpha}}(\gamma,\xi)=\rho_{\xi}^{\check{\alpha}}(\gamma^{-1}o,o)=-\int_o^{\gamma o}\alpha$, whence:
\begin{eqnarray}\label{qgammaforme}
\frac{d(\gamma^{-1}_*\chap{\mu}_o^{\alpha})}{d\chap{\mu}_o^{\alpha}}(\xi,s)=\exp(-\int_o^{\gamma o}\alpha).
\end{eqnarray}
 
\noindent
Remark that $\gamma^{-1}_* \chap{\mu}_o^{\alpha}$ 
is a multiple of $\chap{\mu}_o^{\alpha}$ and does not depend on $o$: 
the cocycle does not depend on $\xi$ and $o$, but only on $\gamma$, and 
the measure $\nu_o^{\alpha}$ on $\D$ is independent on $o$ too. So we will omit the subscript  $o$ 
for these measures $\chap{\mu}^{\alpha}$ et $\nu^{\alpha}$.\\

\noindent
Let us recall some notions: 
the first singular homology group $H_1(M,\Z)$ over $\Z$ can be identified to
the abelianization
$\Gamma/[\Gamma,\Gamma]$ of the fundamental group. 
The first real homology group  
$H_1(M,\R)$ is the $\R$-vector space generated by $H_1(M,\Z)$ except the torsion elements: $H_1(M,\R)=H_1(M,\Z)\otimes_{\Z}\R$.
The first de Rham cohomology group $H^1_{dR}(M)$ is the set of closed $1$-forms modulo exact forms.
Moreover, the de Rham isomorphism (see \cite{warner}) allows to identify the dual of $H_1(M,\R)$ with $H^1_{dR}(M)$. 
This identification is induced by the map:
\begin{eqnarray*} 
  H_1(M,\R)\times H^1_{dR}(M) &\to& \R \\
  (\gamma,\alpha) &\mapsto & \int_{\gamma}\alpha\,.
\end{eqnarray*}
 
\noindent
We denote by $\int_{\gamma}\alpha$ the integral of any representant of 
$\alpha$ on any loop of $M$ belonging to the homology class of $\gamma$.  
In particular, 
it is the integral of $\alpha$ on the unique closed geodesic in the homotopy class of $\gamma$, 
and it is also the integral $\int_o^{\gamma o}\alpha$, since by definition of the action of
$\Gamma$ on $\M$, the curve $[o,\gamma o]$ projects on a loop in the class $\gamma$. 

\noindent
Theorem \ref{1bis} bis allows to determine the measures on ${\cal H}$, which are invariant by a normal subgroup
$\overline{\Gamma}$ of $\Gamma$, and quasi-invariant by $\Gamma$. 
In other words, if $\overline{M}$ is a regular cover of $M$, with fundamental group $\overline{\Gamma}$, 
we can determine the  invariant measures on the space of horospheres of $\overline{M}$, 
which are quasi-invariant by the group $\Gamma/\overline{\Gamma}$ of the cover.
If $\alpha$ is a 
closed $1$-form which vanishes on the image of $\overline{\Gamma}$ in $H_1(M,\R)$, 
formula (\ref{qgammaforme}) above  shows that the measure $\chap{\mu}^{\alpha}$ induced by the Gibbs measure $m^{\alpha}$ is
in particular invariant by $\overline{\Gamma}$. 
Theorem \ref{cover} establishes the converse statement under the additional assumption: if $\chap{\nu}$ is 
$\overline{\Gamma}$-invariant, $\overline{\Gamma}$-ergodic, 
and quasi-invariant by $\Gamma$, then it is necessarily of the form $\chap{\nu}=\chap{\mu}^{\alpha}$, with $\alpha$ a $1$-form on $M$
vanishing on loops of $\bar{\Gamma}$.
%

\begin{demon}{of Theorem \ref{cover}}The first part of the Proposition was proved above. 
The group $\overline{\Gamma}$ is normal in $\Gamma$, so the action of an element
$a\in A=\Gamma/\overline{\Gamma}$ on $\chap{\nu}$ makes sense: indeed, we can write
$a=\gamma_a\overline{\Gamma}=\overline{\Gamma}\gamma_a$, and $\overline{\Gamma}$-invariance of $\chap{\nu}$ shows
that if $\gamma_a\in\Gamma$, the measure ${\gamma_a^{-1}}_*\chap{\nu}$  only depends on $a\in A$.
It shows also that for all $a\in A$, 
the measure $a^{-1}_*\chap{\nu}$ is also $\overline{\Gamma}$-invariant.
The derivative 
$\frac{da^{-1}_*\chap{\nu}}{d\chap{\nu}}$ is then a 
$\overline{\Gamma}$-invariant map $\chap{\nu}$-a.e,  and by ergodicity of $\chap{\nu}$,
it is constant. So for all $a\in A$, there exists a constant
$\lambda(a)>0$, such that $a^{-1}_*\chap{\nu}=\lambda(a)\chap{\nu}$. 
The map $\lambda$ is clearly a  morphism of groups
from $A$ into $\R^*_+$, and when composed with the quotient morphism $\Gamma \to A$, it defines 
a homomorphism, still denoted by $\lambda$, from $\Gamma$ to $\R^*_+$. Moreover, as $\R^*_+$ is abelian, 
$\lambda$  induces a  morphism defined on the abelianization $\Gamma/[\Gamma,\Gamma]=H_1(M,\Z)$.
We can then extend it by $\R$-linearity into the exponential of a linear form on 
$H_1(M,\R)$. 
The de Rham isomorphism between $H_1(M,\R)$ and $H^1_{dR}(M)$ 
allows to see $\log\lambda$ as an element of $H^1_{dR}(M)$; so let $-\alpha$ be a closed $1$-form in its class; 
for all loop $c(\gamma)$ in the homology class of $\gamma$, $\lambda(\gamma)=\exp(-\int_{c(\gamma)}\alpha)$. 
In particular,
$\lambda(\gamma)=\exp(-\int_{\pi[o,\gamma o]}\alpha)=\exp(-\int_o^{\gamma o}\alpha)$. 
Moreover, a closed $1$-form $\alpha$ is in particular a H\"older map on the unit tangent bundle $\tm$, 
hence we have: for all $\gamma\in\Gamma,$ $\gamma^{-1}_*\chap{\nu}=\exp(-\int_o^{\gamma o}\alpha)\,\chap{\nu}$. 
Then Theorem \ref{1bis} bis applies, and $\chap{\nu}$ is necessarily equal to $\chap{\mu}^{\alpha}$.
\end{demon}

\addcontentsline{toc}{section}{Bibliography}



\end{document}